\documentclass[12pt]{article}
\usepackage{amssymb}
\newtheorem{theorem}{Theorem}

\newtheorem{lemma}{Lemma}[section]

\newtheorem{prop}{Proposition}[section]

\newcommand{\qed}{\hfill $\square$\vskip .2cm}

\newcommand{\sect}[1]{\section{#1}\setcounter{equation}{0}}
\newcommand{\comment}[1]{}
\newcommand{\Real}{\mbox{\rm I\kern-.23em\hbox{R}}}

\newcommand{\be}{\begin{equation}}
\newcommand{\ee}{\end{equation}}

\def\h{\mathfrak{h}}

\def\f{\mathfrak{f}}

\def\D{{\cal D}}
\def\L{{\mathfrak{L}}}
\def\l{{\cal L}}

\def\Z{\tilde{Z_\alpha}}

\def\R{\mathbb R}
\def\N{\mathbb N}
\def\Z{\mathbb Z}
\def\S{{\cal S}}
\def\A{{\cal A}}
\def\a{\nabla a}

\newcommand{\<}{\langle}
\renewcommand{\>}{\rangle}
\begin{document}

 \title{On diffusivity of a tagged particle in asymmetric
 zero-range dynamics}
 \author{
\begin{tabular}[l]{ccc}
Sunder Sethuraman\\
Iowa State University
\end{tabular}
}

 \thispagestyle{empty}
 \maketitle
 \abstract{
Consider a distinguished, or tagged particle in zero-range dynamics on
 $Z^d$ with rate $g$ whose
 finite-range jump probabilities $p$ possess a drift $\sum_j jp(j)\neq 0$.  We
 show, in equilibrium, that the variance of the tagged particle
 position at time $t$ is at least order $t$ in all 
$d\geq 1$ and at most
 order $t$ in $d=1$ and $d\geq 3$ for a wide class of rates $g$.
 Also, in $d=1$, when the jump distribution $p$ is totally asymmetric and
 nearest-neighbor, and also when the rate $g(k)$ increases and
 $g(k)/k$ decreases with $k$, we show the diffusively scaled centered
 tagged particle position converges to a Brownian motion.
 }
 
 \vskip .2cm
 {\sl Abbreviated title}:  On diffusivity of a tagged particle in zero-range. \\[.15cm]
 {\sl AMS (2000) subject classifications}: Primary 60K35; secondary 60F05.

\sect{Introduction and Results}
The zero-range process, introduced by Spitzer \cite{Spitzer}, follows
the evolution of a collection of interacting random walks on 
$\Z^d$--namely, from a vertex with $k$ particles, one of the particles
displaces by $j$ with rate $g(k)p(j)$.  The function on the
non-negative integers $g:\N \rightarrow \R_+$ is called the process ``rate,'' and $p(\cdot)$ denotes the translation-invariant single
particle transition probability.
The above interaction is in the ``time-domain'' but not ``spatially,'' 
hence the name ``zero-range.''  We note the case when $g(k) \equiv k$
describes the situation of completely independent particles.

More precisely, let $\Sigma = \N^{\Z^d}$ be the configuration space
where a configuration $\xi = \{\xi_i:i\in Z^d\}$ is given through
occupation numbers $\xi_i$ at vertex $i$.  The zero-range system then
is a Markov process $\xi(t)$ on the space of right-continuous paths
with left limits $D(\R_+,\Sigma)$ with formal generator defined on
real test functions $\phi$,
$$(L\phi)(\xi) \ = \ \sum_j\sum_i g(\xi_i)p(j-i)
(\phi(\xi^{i,j})-\phi(\xi))$$
where $\xi^{i,j}$ is the configuration where a particle from $i$ is
moved to $j$.  That is, $\xi^{i,j} = \xi -\delta_i +\delta_j$ where
$\delta_k$ is the configuration with a single particle at $k$.

When a particle is distinguished, or tagged, we can consider the joint
process $(x(t),\xi(t))$ on $D(\R_+, \Z^d\times \Sigma)$ where $x(t)$ is the position of the tagged
particle at time $t$.  The formal generator is given by
\begin{eqnarray*}
(\L\psi)(x,\xi) &=& \sum_j\sum_{i\neq x} g(\xi_i)p(j-i)
(\psi(x,\xi^{i,j})-\psi(x,\xi))\\
&&\ \ + \sum_j g(\chi_x)\frac{\xi_x -
  1}{\xi_x}p(j)(\psi(x,\xi^{x,x+j})-\psi(x,\xi))\\
&&\ \ + \sum_j
\frac{g(\xi_x)}{\xi_x}p(j)(\psi(x+j,\xi^{x,x+j})-\psi(x,\xi)).
\end{eqnarray*}
Here, the first term corresponds to particles other than at the tagged
particle position $x$ moving, the second term corresponds to other
particles moving from $x$, and the last term represents motion of the
tagged particle itself.

It will be convenient to consider the ``reference'' process from the
point-of-view of the tagged particle, that is $\eta(t) =
\{\xi_{i+x(t)}(t): i\in \Z^d\}$ which can be
obtained from the map $\pi((x(\cdot),\xi(\cdot)))= \eta(\cdot)$,
and
has formal generator
\begin{eqnarray*}
(\l\phi)(\eta) &=&  \sum_j\sum_{i\neq 0} g(\eta_i)p(j-i)
(\phi(\eta^{i,j})-\phi(\eta))\\
&&\ \ + \sum_j g(\eta_0)\frac{\eta_0 -
  1}{\eta_0}p(j)(\phi(\eta^{0,j})-\phi(\eta))\\
&&\ \ + \sum_j
\frac{g(\eta_0)}{\eta_0}p(j)(\phi(\tau_j(\eta^{0,j}))-\phi(\eta)).
\end{eqnarray*}
Here, $\tau_j(\eta^{0,j})$ is the configuration obtained by displacing the
tagged particle by $j$ and then shifting the reference frame to its
position; the notation $(\tau_j\eta)_k = \eta_{k+j}$ for $k\in \Z^d$.

The construction of these systems requires some conditions on $g$ and
$p$.
Namely, we will assume throughout $g(0)=0$, $g(k)>0$ for
$k\geq 1$, $|g(k+1)-g(k)|\leq K$ for some constant $K$, and
$\liminf_{k\rightarrow \infty} g(k)>0$, and also $p$ is
finite-range, that is $p(i) = 0$ for $|i|\geq R$ for some $1\leq R<\infty$,
whose symmetrization $s(x)=(p(x)+p(-x))/2$ is irreducible.
Under weaker assumptions, which include the above, Andjel constructs
the process $\xi(t)$ semigroup $T^L_t$ on a class of ``Lipschitz''
functions $\D$ defined on a subset $\Sigma' \subset \Sigma$ of the
configuration
space,
$$\Sigma' \ = \ \bigg\{\xi: \|\xi\| = \sum_{i\in \Z^d} |\xi_i|\beta_i <\infty\bigg\}$$
$$\D \ = \ \bigg\{f: |f(\xi')-f(\xi'')| \leq c\|\xi' - \xi''\| \ {\rm
  for \ all \ }\xi',\xi''\in \Sigma', \ {\rm for \ some \
  }c=c(f)\bigg\}$$
where one can take $\beta_i = \sum_{n\geq 0} 2^{-n}s^{(n)}(i)$ for
  instance \cite{Andjel}.
In a similar way, one can
  construct the process $(x(t),\xi(t))$ semigroup $T^\L_t$ 
with respect to ``Lipschitz'' functions $f$ where
  $|f(x,\xi')-f(y,\xi'')|\leq c[|x-y| +\|\xi'-\xi''\|]$ for all
  $x,y\in \Z^d$ and $\xi',\xi''\in \Sigma'$.  Then, also from the
  map $\pi$, process
  $\eta(t)$ semigroup $T^\l_t$ can be constructed on $\D$.  

The zero-range process $\xi(t)$ has a well known explicit family
product invariant
measures $R_\alpha = \prod_{i\in Z^d} \mu_\alpha$ for
$0\leq\alpha<\liminf g(k)$ with marginal
$$\mu_\alpha(k) \ = \ \frac{1}{Z_\alpha}\frac{\alpha^k}{g(1)\cdots g(k)}
 {\rm \ \ when \ \ } k\geq 1 \ \ \ {\rm and \ }
\mu_\alpha(0) \ = \ \frac{1}{Z_\alpha}  \ {\rm \ when \ \ }k=0.
$$
where $Z_\alpha$ is the normalization \cite{Andjel}. 
Let $\rho(\alpha) = \sum_k k\mu_\alpha(k)$ be the density of particles
under $R_\alpha$, and let $\rho^*= \lim_{\alpha \uparrow \liminf
  g(k)}$ [note that $\rho^*$ may be finite for some type of $g$'s].
As $\rho(\alpha) \uparrow \rho^*$ for 
$\alpha\uparrow \liminf g(k)$, for a given $0\leq \rho<\rho^*$, 
there is a unique
inverse $\alpha = \alpha(\rho)$.

For the reference process $\eta(t)$, the ``palm'' measures
given by $dQ_\alpha = (\eta_0/\rho(\alpha))dR_\alpha$ are 
invariant (cf. \cite{Port}, \cite{Saada}).  Only the marginal at the
origin, which we denote $\mu^0_\alpha$, differs
$$\mu^0_\alpha(k) \ = \ \frac{1}{Z_\alpha}\frac{k}{\rho}\frac{\alpha^{k}}{g(1)\cdots
  g(k)}
\ \ {\rm for \ \ }k\geq 1.$$

We now remark that, with respect to an invariant $R_\alpha$, 
one can extend the zero-range process semigroup
$T^L_t$ and generator $L$ to $L^2(R_\alpha)$ so that bounded functions
in $\D$ form a core (cf. section 2 \cite{Sext}).  In the same way,
with respect to a $Q_\alpha$, the reference semigroup $T^\l_t$ and
generator $\l$ can be extended to $L^2(Q_\alpha)$ with the same core.
We note also here constructions of these processes can be made through
the martingale-problem approach \cite{Saada}, \cite{Saada-thesis}.
Also, in this context, we note a Hille-Yosida type approach \cite{Liggett1}.

In addition, we note
both families $\{R_\alpha\}$ and $\{Q_\alpha\}$ are in fact extremal
measures in their respective convex set of invariant measures, and so 
process evolutions starting from such invariant states are
time-ergodic \cite{Sext}.  Also, we note the adjoints $L^*$ and $\l^*$
with respect to $R_\alpha$ and $Q_\alpha$ respectively correspond to
``time-reversal'' and are themselves zero-range and
reference processes but with reversed jump probabilities
$p^*(\cdot )=p(- \cdot)$.  Finally, for $0\leq \alpha <\liminf g(k)$,
we note $\mu_\alpha$ and $\mu^0_\alpha$ possess all moments.

In the following, to avoid degeneracies, we will work with a fixed
$0<\alpha<\liminf g(k)$ for which $\rho(\alpha)>0$, and
corresponding $R_\alpha$ and $Q_\alpha$.  
For simplicity, we will denote by $E_\alpha [\cdot]$ and
$P_\alpha[ \cdot]$ the expectation and probability for the process measures
starting from $Q_\alpha$ when there is no confusion; otherwise, the
underlying measure will noted as a suffix.

We now discuss the problem studied in this article and its history.
The question of tagged
particle asymptotics
was even mentioned in Spitzer's seminal paper.  Such questions are
important to physics and other applications \cite{LebSpohn}.
What is known are
some laws of
large numbers (LLN), and some equilibrium central limit theorems (CLT)
in ``local balance'' cases.  

Write, with respect to the reference process, the position $x(t)$ as the sum
total displacement ``shift,'' that is
$$x(t) \ = \ \sum_j jN_j(t)$$
where $N_j(t)$ is the number of ``shifts'' of size $j$ the reference
process makes up to time $t$.
The count
$N_j(t)$ is compensated by $\int_0^t (g(\eta_0(s)/\eta_0(s))p(j)ds$,
so that further
\begin{equation}
\label{martdecomp}
x(t) \ = \ \sum_j jM_j(t) 
+\sum_j j\int_0^t \frac{g(\eta_0(s))}{\eta_0(s)}p(j)ds
\end{equation}
where $M_j(t) = N_j(t) -\int_0^t (g(\eta_0(s)/\eta_0(s))p(j)ds$ are
orthogonal martingales.  Moreover, we note $M_j^2(t) - \int_0^t
(g(\eta_0(s))/\eta_0(s)) p(j)ds$ are also martingales.

So, the tagged position $x(t)$ is a function
of the reference process.  For most of the paper, we will use this
``reference frame'' interpretation, that is, the notation $x(t)$ will
denote the ``compound shift'' $\sum jN_j(t)$.
It will also be useful to define
$$M(t) \ = \ \sum jp(j)M_j(t) \ \ \ {\rm and \ \ \ }A(t) \ = \ 
\int_0^t \f(\eta(s))ds$$
where $\f(\eta) = (\sum jp(j))(g(\eta_0)/\eta_0) - (\alpha/\rho)$.

Then, in equilibrium, that is when the reference process is under initial distribution
$Q_\alpha$, one obtains 
$$E_\alpha[x(t)] \ = \ E_\alpha[g(\eta_0)/\eta_0]\sum jp(j) \ = \
\frac{\alpha}{\rho(\alpha)}\sum jp(j)$$
and LLN
$$\lim_{t\rightarrow \infty} \frac{1}{t} x(t) \ = \
\frac{\alpha}{\rho(\alpha)} \sum_{j} jp(j) \ \ \ {\rm a.s.  }$$
(cf. \cite{Saada}, \cite{Sext}).
Also, we refer the reader to some interesting
LLN results under some non-equilibrium initial distributions
\cite{Reza-lln}.  

With respect to fluctuations, when the jump probabilities are mean-zero,
$\sum jp(j)=0$, then $x(t) = \sum jM_j(t)$ is a martingale
as the compensator terms cancel.  Under equilibrium $Q_\alpha$, 
the quadratic variation is 
$$E_\alpha[|x(t)|^2] \ = \ \sum |j|^2 \int_0^t \frac{g(\eta_0(s))}{\eta_0(s)}
p(j) ds \ \rightarrow \ \frac{\alpha}{\rho}\sum j^2 p(j) \ \ \ {\rm
  a.s.}$$
and so by martingale central limit theorem one gets the invariance
principle 
$$\lim_{\lambda \rightarrow \infty} \frac{1}{\sqrt{\lambda}}x(\lambda
t) \ = \ B_\alpha(t)$$
where $B_\alpha(t)$ is $d$ dimensional Brownian motion with covariance
matrix $((\alpha/\rho)t\sum_j (e_i\cdot j)(e_k \cdot j)p(j))$ where
$\{e_i\}$ is the standard basis of $Z^d$ \cite{Saada}, \cite{Sext}.

The goal of this article is to further characterize the equilibrium
fluctuations when the jump probability has a drift $\sum jp(j)\neq
0$.  The first result is that the tagged particle variance is at least
diffusive in all dimensions without conditions.  As a comparison, 
we note this is not
true for simple exclusion in the case $d=1$ and the jump probability
$p$ is nearest-neighbor symmetric where the variance at time $t$ is
order $t^{1/2}$ \cite{Arratia}.
\begin{theorem}
Under initial distribution $Q_\alpha$, we have in all dimensions
$d\geq 1$ for $t\geq 0$ that
$$\bigg[\frac{\alpha}{\rho(\alpha)}\sum_j |j|^2 p(j)\bigg ] t 
\ \leq \ E_\alpha\bigg[
|x(t) - 
E_\alpha[x(t)]|^2\bigg ].
$$
\end{theorem}

{\it Proof.} These bounds
follow from an explicit calculation.
We have
\begin{eqnarray}
\label{martsequence}
E_\alpha[|x(t)-E_\alpha[x(t)]|^2]
&=& E_\alpha [|M(t) + A(t)|^2]\nonumber \\
&=& E_\alpha [|M(t)|^2] + 2E_\alpha [M(t)\cdot A(t)] + E_\alpha
[|A(t)|^2] \\
&=& \frac{\alpha}{\rho}\sum |j|^2p(j)t + 2\int_0^t E_\alpha [M(s)\cdot \f(\eta(s))]ds
+ E_\alpha [|A(t)|^2].\nonumber
\end{eqnarray}
Now, under time reversal at $s$, 
$\eta^*(u)= \eta(s-u)$,
the number of $j$-shifts up to time $s$ equals the number of
$-j$-shifts in the reversed process up to time $s$, 
$N_j(s; \eta(\cdot)) = N_{-j}(s; \eta^*(\cdot))$.  Also, 
$M^*(s) = \sum j N_{-j}(s) - \sum j\int_0^s
(g(\eta^*_0(u))/\eta^*_0(u)) p(j)ds$ is a martingale with respect to
the reversed process.  So,
we have
$E_\alpha [M(s)\cdot\f(\eta(s))] =
E_\alpha [ M^*(s) \cdot \f(\eta^*(0))] = 0$.
Hence,
\begin{equation}
\label{v_decomp}
E_\alpha[|x(t)-E_\alpha[x(t)]|^2] \ = \ \frac{\alpha}{\rho}\sum
|j|^2p(j)t
 + E_\alpha [|A(t)|^2] \ \geq \ \frac{\alpha}{\rho}\sum
|j|^2p(j)t
\end{equation}
(we remark a different representation holds for exclusion processes \cite{DeMF}).
\qed

To give some upperbounds on the tagged particle variance, we describe
some classes of rate functions $g$.

{\bf Assumption (SP).} Let $L_n$ be the generator of the symmetric zero-range
process on a cube $B_n=\{i\in \Z^d: |i|\leq n\}$, namely
$(L_n\phi)(\xi) 
= \sum_{i,j\in B_n}g(\xi_i)(\phi(\eta^{i,j})-\phi(\eta))s(j-i)$.  Let
$W(n,M)$ be the inverse of the spectral gap of $L_n$ when there are
$M$ particles in $B_n$.  
Then, we assume the rate function $g$ is such
that there is a constant $C=C(\alpha,p,d)$ where
$E_{R_\alpha}[(W(n,\sum_{i\in B_n}\xi_i))^2] \leq Cn^4$.
\vskip .1cm

We observe rates $g$ where $W(n,M)\leq Cn^2$ for a constant $C=C(d)$
independent of $M$, satisfy (SP) trivially, and include those rates where, for
some $a\geq 1$ and $b>0$, $g(k+a)- g(k)\geq b$ for all $k\geq 0$ \cite{LSV}.
Also, for the rate $g(k)=1_{[k\geq 1]}$, it is known $W(n,M) \leq
C(1+M/n)^2 n^2$ for some constant $C=C(d)$ \cite{Morris}, and so
(SP) holds.  It is most likely true that all rates $g$ satisfy
(SP).

\vskip .2cm
{\bf Assumption (ID).}  The rate function $g$ is such that $g(k)$
increases and $g(k)/k$ decreases with $k$.
\vskip .1cm



\begin{theorem}
\label{thm1}
Under initial distribution $Q_\alpha$, when $\sum_i ip(i)\neq 0$, we have in $d=1$ (without
special assumptions), and in $d\geq 3$ under Assumption (SP)
that there is a constant $C=C(\alpha,p,d)$ where for $t\geq 0$
$$ E_\alpha\bigg[
|x(t) - 
E_\alpha[x(t)]|^2\bigg ] \ \leq Ct.$$
\end{theorem}

We also note an invariance principle in a special case in $d=1$.
\begin{theorem}
\label{thm2}
Under initial distribution $Q_\alpha$, in $d=1$ when the jump
probability is totally asymmetric $p(1)=1$ and $g$ satisfies
Assumption (ID), we have
the invariance principle
$$\lim_{\lambda \rightarrow \infty} \frac{1}{\sqrt{\lambda}} \bigg
(x(t) - E_\alpha[x(t)]\bigg) \ = \ B_\alpha(t)$$
where $B_\alpha$ is Brownian motion with diffusion coefficient
$\sigma^2(\alpha)> \alpha/\rho$.
\end{theorem}

To compare, we remark with respect to asymmetric simple exclusion
similar invariance principles have been shown
in $d\geq 3$ for finite-range $p$ \cite{SVY} and in $d=1$ when
$p$ is nearest-neighbor \cite{Kipnis}.
In this context, perhaps the main contribution of this paper is
Theorem 
\ref{thm1} in $d=1$ as these upperbounds, valid for the general
finite-range zero-range process, have no counterpart in the simple
exclusion results.
 
The proof of
Theorem \ref{thm1} follows from an analysis of certain variance or
$H_{-1}$ norms, and is found in section 2.  The proof of
Theorem \ref{thm2}, in section 3, shows that a tagged particle has positively
correlated increments in the totally asymmetric nearest-neighbor case
in $d=1$.  Combined with diffusive variance bounds (Theorem
\ref{thm1}), the invariance principle follows by applying a
Newman-Wright theorem.  We note, in comparison, a tagged particle in  
$d=1$ simple exclusion
with totally asymmetric nearest-neighbor transitions has negatively
correlated increments, and in fact is a Poisson process
\cite{Liggett}.  In the zero-range context, since $\sigma^2(\rho)
>\alpha/\rho = E_\alpha[x(1)]$, the tagged particle is not a Poisson
process.

\sect{Proof of Theorem \ref{thm1}}

We first discuss some definitions and estimates involving variational
formulas for some resolvent quantities.
Note that $\D\subset L^2(Q_\alpha)$ as for $f\in \D$ we have
\begin{eqnarray*}
E_\alpha[|f(\eta)|^2] &\leq & c^2 E_\alpha\bigg[\big(\sum \eta_i\beta_i\big)^2\bigg]\\
& \leq & c^2\sum_{i\neq j} E_\alpha[\eta_i \eta_j]\beta_i\beta_j + c^2\sum
E_\alpha[\eta_i^2] \beta_i^2\\
&\leq& c^2 C\bigg (\sum \beta_i \bigg)^2 \ < \ \infty
\end{eqnarray*}
for a constant $C=C(\alpha)$.  Also, we note the notation $E_\alpha[fg] = \<f,g\>_\alpha$, and
$E_\alpha[f^2] = \|f\|^2_0$.

The generator $\l$ can be decomposed into symmetric and anti-symmetric
parts,
$\l=\S+\A$ where $\S=(\l +\l^*)/2$ and $\A=(\l-\l^*)/2$.  Consider the
resolvent operator $(\lambda -\l)^{-1}:L^2(Q_\alpha) \rightarrow
L^2(Q_\alpha)$ 
well defined for $\lambda>0$; in particular, $(\lambda -\l)^{-1} f =
\int_0^\infty e^{-\lambda s}(T^{\l}_s f) ds$.  Since the symmetrization of $(\lambda
-\l)^{-1}$ has inverse $(\lambda - \l^*)(\lambda -\S)^{-1}(\lambda -\l) =
(\lambda -\S) +\A^*(\lambda-\S)^{-1}\A$, we have
the variational formula for $f\in \D$,
$$\<f,(\lambda-\l)^{-1} f\>_\rho \ = \ \sup_{\phi\in \D} \bigg\{
2\<f,\phi\>_\alpha - \<\phi, (\lambda -\S)\phi\>_\alpha -\<\A\phi, (\lambda
- \l)^{-1}\A\phi\>_\alpha\bigg \}.$$
Now, as $\A^* = -\A$ and $\A^*(\lambda -\S)^{-1}\A$ is a non-positive
operator,
we have the easy bound that $\<f,(\lambda-\l)^{-1}f\>_\alpha$ is bounded
by its ``symmetrization,''
\begin{eqnarray}
\<f,(\lambda-\l)^{-1} f\>_\alpha & \leq & \sup_{\phi \in \D} \bigg\{
2\<f,\phi\>_\alpha - \<\phi, (\lambda -\S)\phi\>_\alpha \bigg\}\nonumber\\
&=& \<f,(\lambda -\S)^{-1}\>_\alpha.
\label{symbound}
\end{eqnarray}

It will be convenient to
define, for $f\in\D$, the $H_1(Q_\alpha)$ (semi)-norm by $\|f\|^2_1 =
\<f, (-\l)f\>_\alpha = \<f, (-\S)f\>_\alpha$.  The $H_1$ space then is the
completion with respect to this norm.  Explicitly, for $\psi\in \D$,
\begin{eqnarray*}
\<\psi, (-\S)\psi\>_\alpha &=& \frac{1}{2}\sum_j\sum_{i\neq 0}
E_\alpha[g(\eta_i)(\psi(\eta^{i,i+j})-\psi(\eta))^2] s(j)\\
&&\ \ \ \ +\frac{1}{2}\sum_j E_\alpha[g(\eta_0)\frac{\eta_0
  -1}{\eta_0}((\psi(\eta^{0,j})-\psi(\eta))^2] s(j)\\
&&\ \ \ \  + \frac{1}{2}\sum_j
E_\alpha[\frac{g(\eta_0)}{\eta_0}(\psi(\tau_j(\eta^{0,j})) -
\psi(\eta))^2]s(j).
\end{eqnarray*}
Let $H_{-1}$ be the dual of $H_1$, namely, the completion over $\D$ with respect to norms
$\|\cdot\|_{-1}$ given in terms of variational formulas
\begin{equation}
\|f\|^2_{-1} \ = \ \sup_{g\in \D} \bigg\{2\<f,g\>_\alpha -
\|g\|_1^2\bigg\}
\ =\  \sup_{g\in\D} \frac{\<f,g\>_\alpha }{ \<g,(-\S) g\>^{1/2}_\alpha}.
\label{var_expression}
\end{equation}
Similarly, define the notation
$\|f\|^2_{1,\lambda} =\<f,(\lambda-\S)f\>_\alpha$ and
$\|f\|^2_{-1,\lambda} =\sup_{g\in\D} \{2\<f,g\>_\alpha -
\|g\|^2_{1,\lambda}\}=\<f,(\lambda -\S)^{-1}f\>_\alpha$.
Note also, clearly as one can drop the $\lambda\<f,f\>_\alpha$ term, 
that $\|f\|^2_{-1,\lambda} \leq \|f\|_{-1}^2$.
In addition, we note the following useful ``resolvent'' estimate.  For
$f\in L^2(Q_\alpha)$ we have
\begin{equation}
\label{resolventbound}
\|f\|^2_{-1,\lambda} \ =\  \int_0^\infty e^{-\lambda
  t}\<T^{\S}_tf,f\>_\alpha dt \ \leq\  \frac{1}{\lambda} \|f\|^2_0
\end{equation}
where $T^\S_t$ is the semigroup for the symmetrized process.
Then, for $g\in \D$, we have
$$E_\alpha[fg] \ \leq \ \|f\|_{-1,\lambda}\|g\|_{1,\lambda}
\ \leq\ \bigg[\frac{1}{\lambda}\|f\|^2_0\bigg]^{1/2} \|g\|_{1,\lambda}.
$$

\vskip .2cm
Now, for $f\in L^2(Q_\alpha)$, let 
$\sigma^2_t(f) = E_\alpha[(\int_0^t f(\eta(s))ds)^2]$, and observe
from the decomposition (\ref{v_decomp}), to get diffusive bounds on the tagged particle
variance, one need only bound 
$$\sigma^2_t(\h) \ = \ E_\alpha \bigg[ \bigg(\int_0^t
\h(\eta(s))
ds\bigg)^2\bigg]
\ <\ Ct
$$
where $\h(\eta) = (g(\eta_0)/\eta_0) - (\alpha/\rho)$ for some
constant $C=C(\alpha,p,d)$.
The next result relates $\sigma^2_t(f)$ to some $H_{-1}$ norms.

\begin{prop} 
\label{ub1} For $f\in L^2(Q_\alpha)$, there is a universal constant $C_1$ such that 
\begin{eqnarray*}
\sigma^2_t(f) &\leq& C_1t\<f, (t^{-1} - \l)^{-1}
f\>_\alpha\\
&\leq&C_1t\<f, (t^{-1} - \S)^{-1} f\>_\alpha \ \leq \ C_1 t \|f\|_{-1}^2.
\end{eqnarray*}
\end{prop}

{\it Proof.} The first line is well-known (with a proof found for
instance in Lemma 3.9 \cite{Sclt}), the second bound is
(\ref{symbound}), and the third bound is explained after (\ref{var_expression}). \qed

\vskip .2cm
{\it Proof of Theorem \ref{thm1}.}
The strategy to bound $\sigma^2_t(\h)$ falls into
two cases $d=1$ and $d\geq 3$ under (SP).
We first comment on the case $d=1$, and then on the $d\geq 3$ case.

\vskip .1cm
{\it Case $d=1$.}  (1) We
will find a sequence of functions (in subsection 2.1.1)
$\{\phi_\lambda: 0<\lambda \leq 1\}\subset\D$ such that 
\begin{equation}
\label{bounds1}
\sup_{0<\lambda \leq 1} \|\h -\l\phi_\lambda\|_{-1,\lambda} \ < \
\infty
\end{equation}
and also
\begin{equation}
\label{bounds2}
\sup_{0<\lambda \leq 1} \bigg (\|\phi_\lambda\|^2_1 +
\lambda\|\phi_\lambda\|^2_0\bigg ) \ < \infty.
\end{equation}

(2) Note $M_t(f) = f(\eta(t))-f(\eta(0)) - \int_0^t (\l f)(\eta(s))ds$
is a martingale for $f\in\D$ with quadratic variation (by stationarity)
$E_\alpha [ (M_t(f))^2] = 2tE_\alpha[ f(-\l)f] = 2t\|f\|_1^2$.  Then,
we can write
$$-\int_0^t \l\phi_\lambda (\eta(s))ds \ = \ M_t(\phi_\lambda) +
\phi_\lambda (\eta(0))-\phi_\lambda (\eta(t))$$
and so (by stationarity) 
$$
\sigma^2_t(\phi_\lambda) \ \leq \  
6\bigg(t\|\phi_\lambda\|^2_1 +  \|\phi_\lambda\|^2_0\bigg )
\ = \ 6t\bigg (\|\phi_\lambda\|^2_1 +
\frac{1}{t}\|\phi_\lambda\|^2_0\bigg ).
$$

(3) Hence, by choosing $\lambda = t^{-1}$, we have from Proposition
    \ref{ub1} that
\begin{eqnarray}
\sigma^2_t(\h) &\leq& 2\sigma^2_t(\h -\l\phi_{t^{-1}})
+2\sigma^2_t(\l\phi_{t^{-1}})\nonumber\\
&\leq& 2C_1 t \|\h - \l\phi_{t^{-1}}\|^2_{-1,t^{-1}} + 12t\bigg(\|\phi_{t^{-1}}\|^2_1 +
\frac{1}{t}\|\phi_{t^{-1}}\|^2_0\bigg ).
\label{sequence}
\end{eqnarray}
Then, by (1) and (2), $\sigma^2_t(\h) \leq Ct$ for some constant
$C=C(\alpha,p)$ and $t\geq 1$.  
For $0\leq t<1$, bounds are immediate.
This finishes the proof in
this case.

\vskip .1cm
{\it Case $d\geq 3$ and (SP).}
By Lemma \ref{ub1}, we need only show $\|\h\|_{-1}<\infty$.  One may
be able to do this directly by ``integration-by-parts'' but as the
$Q_\alpha$ marginal at the origin differs from the other marginals,
one cannot apply immediately results in the literature.  So, we ``modify''
the function $\h$ and then apply these results.  

Let $j_0$ be a point in the support of $p$.  Consider the function
$\phi(\eta)=(\eta_{j_0}-\rho)/(\rho p(j_0)) \in\D$.  In subsection 2.1.2, we show
that 
$\|\h-\l\phi\|_{-1} <\infty$.  Clearly $\|\phi\|_1<\infty$ and
$\|\phi\|_0<\infty$.  
Then, by following the sequence
(\ref{sequence}), we have
$$
\sigma^2_t(\h) \  \leq \ 2C_1t \|\h-\l\phi\|_{-1} + 12t\bigg(\|\phi\|^2_1 +
\frac{1}{t}\|\phi\|^2_0\bigg )
\ <\ Ct
$$
for a constant $C=C(\alpha,p,d)$ and $t\geq 1$.  Bounds when
$0\leq t<1$ are clear.
This finishes the proof. \qed

\subsection{Some Estimates}
We now turn to supplying the needed estimates in the two cases.  We first make a
calculation valid in any dimension $d\geq 1$.
Let
$$\phi(\eta) \ = \ \sum_{i\in Z^d} a_i (\eta_i-\rho)$$
where $\sum a_i^2 <\infty$.  Clearly $\phi\in L^2(Q_\alpha)$, and
$\phi=\lim \phi^n$ is the $L^2(Q_\alpha)$ limit of functions $\phi^n =\sum_{|i|\leq
  n}a_i(\eta_i-\rho)\in \D$.  Also, computes $\l\phi$ as the
$L^2(Q_\alpha)$ limit $\l\phi = \lim \l\phi^n$ as bounded functions in
$\D$ are a core.  
To this end, for $n$ large, observe
$$\phi^n(\eta^{i,i+j}) - \phi^n(\eta) \ = \ a_{i+j} - a_i$$
and, as $\tau_j(\eta^{0,j}) = \tau_j(\eta +\delta_j - \delta_0) =
\tau_j \eta +\delta_0 - \delta_{-j}$, and so $\phi^n(\tau_j
\eta^{0,j}) = \sum_{|i|\leq n}
a_i (\eta_{i+j}-\rho) + a_{0}-a_{-j}$,
we have
$$ \phi^n(\tau_1 \eta^{0,1}) - \phi^n(\eta) \ = \ \sum_{|i|\leq n} (a_i-a_{i+j})
(\eta_{i+j} - \rho) +(a_{0} - a_{-j}).$$
These computations enable us to write
\begin{eqnarray*}
(\l\phi)(\eta) & = & \sum_j\sum_{i\neq 0} (a_{i+j} - a_i)g(\eta_i)p(j)
+ \sum_j(a_{j} - a_0)g(\eta_0)\frac{\eta_0 -1}{\eta_0}p(j) \\
&&\ \ - \sum_j\sum_i (a_{i+j}- a_{i})(\eta_{i+j} - \rho)\frac{g(\eta_0)}{\eta_0}p(j) + \sum_j(a_{0}
-a_{-j})\frac{g(\eta_0)}{\eta_0}p(j)\\
&=& \sum_j\sum_{i\neq 0,-j} (a_{i+j} - a_i)\big [ g(\eta_i) -
{\eta_{i+j}}\frac{g(\eta_0)}{\eta_0}\big ]p(j)\\
&&\ \  + \sum_j(a_0 - a_{-j})\big [
g(\eta_{-j}) - g(\eta_0)\frac{\eta_0 -1}{\eta_0} \big ]p(j)\\
&&\ \  + \sum_j(a_{j} -
a_0)g(\eta_0)\big [ 1 - \frac{\eta_{j}+1}{\eta_0}\big ]p(j).
\end{eqnarray*}
Here, we used $\sum_j\sum_i(a_{i+j}-a_i) = 0$ to reduce the first
sum in the second line.

We now note the
following basic useful computations.
\begin{lemma}
\label{computations}
Let $k\in \Z^d$ be a non-zero vertex, $k\neq 0$.
Then, for $\psi\in L^2(Q_\alpha)$, we have
$$E_\alpha[g(\eta_k)\psi(\eta)] = \alpha E_\alpha[\psi(\eta
+\delta_k)], \  E_\alpha[g(\eta_0)\frac{\eta_0-1}{\eta_0}
\psi(\eta)]=\alpha E_\alpha[\psi(\eta +\delta_0)],$$
and
$$E_\alpha[\frac{g(\eta_0)}{\eta_0}\psi(\eta)] = 
E_\alpha[\frac{g(\eta_0)}{\eta_0}\psi(\tau_j(\eta^{0,j}))], \ 
E_\alpha[(\eta_j +1)\frac{g(\eta_0)}{\eta_0}\psi(\eta)] = E_\alpha[g(\eta_0)\psi(\tau_{-j}(\eta^{0,-j}))].$$
\end{lemma}

{\it Proof.}  We show the last equality as the others are similar.
Write
 \begin{eqnarray*}
E_\alpha [\psi(\eta)(\eta_j +1)\frac{g(\eta_0)}{\eta_0}] &=& \frac{\alpha}{\rho} E_{P_\alpha}[\psi(\eta
 +\delta_0)(\eta_j+1)]
 = \frac{1}{\rho}E_{P_\alpha}[g(\eta_j)\psi(\eta^{j,0})\eta_j]\\
&=&
 E_\alpha[g(\eta_0)\psi((\tau_{-j}\eta)^{j,0}] =
 E_\alpha [g(\eta_0)\psi(\tau_{-j}(\eta^{0,-j}))].
\end{eqnarray*}
\qed

Let now $\psi\in \D$ be a function.  We can write, with Lemma \ref{computations},
\begin{eqnarray}
 E_\alpha[(\l\phi) \psi]
  & =& \sum_j\sum_{i\neq 0,-j} (a_{i+j}-a_i)E_\alpha [(g(\eta_i) -
{\eta_{i+j}}\frac{g(\eta_0)}{\eta_0})\psi(\eta)]p(j) \nonumber \\
 && \ \ + \alpha \sum_j(a_0-a_{-j})E_\alpha [\psi(\eta +\delta_{-j}) - \psi(\eta
 +\delta_{0})]p(j)\nonumber\\
 &&\ \ + \sum_j(a_{j}-a_0)E_\alpha [g(\eta_0)(\psi(\eta ) - \psi(\tau_{-j}( \eta^{0,-j})))]p(j).
\label{lphipsi} \end{eqnarray}
It will be convenent, for later purposes, to observe that in the above
computation we can take $E_\alpha [\psi]=0$ without loss of generality
as $E_\alpha [\l\phi]=0$.

\vskip .2cm
{\bf 2.1.1 Estimates in $d=1$.}  We now work in dimension $d=1$, 
and choose the sequence
$$a_i \ = \ \left\{\begin{array}{rl}
0& \ {\rm for \ }i\leq 0\\
c(1-\lambda)^{i-1}& \ {\rm for \ } i\geq 1
\end{array}\right.
$$
with $c= \rho^{-1}$.
For ease of notation, define $\a_{k,j} = a_{k} - a_j$ and 
note
$$\a_{i+j,i} \ = \ \left\{\begin{array}{rl}
0& \ {\rm for \ } i,i+j\leq 0\\
c(1-\lambda)^{i+j-1}& \ {\rm for \ }i+j\geq 1 \ {\rm and \ }i\leq 0\\
-c(1-\lambda)^{i-1}& \ {\rm for \ }i\geq 1 \ {\rm and \ }i+j \leq 0\\
c[(1-\lambda)^{j}-1] (1-\lambda)^{i-1}& \ {\rm for \ }i,i+j\geq 1.
\end{array}
\right.
$$
Clearly $|\a_{i+j,i}| \leq |c|$
for all $i,j\in \Z$.

Recall now
the range $R$ of the
 distribution $p$, and write, with Lemma \ref{computations},
\begin{eqnarray*}
E_\alpha [(\l\phi)(\eta)\psi(\eta)]
&=& \sum_j\sum_{i\geq R+1} \a_{i+j,i}E_\alpha [(g(\eta_i) -
{\eta_{i+j}}\frac{g(\eta_0)}{\eta_0})\psi(\eta)]p(j)\\
&&\ \ +\sum_j\sum_{{|i| \leq R} \atop {i\neq 0,-j}} \a_{i+j,i}E_\alpha[(g(\eta_i) -
{\eta_{i+j}}\frac{g(\eta_0)}{\eta_0})\psi(\eta)]p(j) \\
 && \ \ +
 \sum_j\a_{0,-j}E_\alpha[g(\eta_0)\frac{\eta_0-1}{\eta_0}(\psi(\eta^{0,-j}) - \psi(\eta))
]p(j)\\
&&\  \ 
 + \sum_j\a_{j,0}E_\alpha[g(\eta_0)(\psi(\eta ) - \psi(\tau_{-j}
 \eta^{0,-j}))]p(j)\\
&=& I_1 +I_2 +I_3 +I_4.
\end{eqnarray*}
Consider now the term $I_1$.  Since $E_\alpha[\psi]=0$, we can write
\begin{eqnarray*}
I_1&=&
\sum_j\sum_{i\geq R+1} \a_{i+j,i}E_\alpha[(g(\eta_i) - \alpha)
-
{(\eta_{i+j}-\rho)}\frac{g(\eta_0)}{\eta_0})\psi(\eta)]p(j)\\
&&\ \ \ \ \ \ \ -\rho\sum_j\sum_{i\geq R+1}
\a_{i+j,i}E_\alpha[\frac{g(\eta_0)}{\eta_0} \psi(\eta)]p(j)\\
&=& J_1 +J_2.
\end{eqnarray*}
Note the last term $J_2$ equals, using $E_\alpha[\psi]=0$ again,
\begin{eqnarray*}
J_2 & = & 
\rho a_{R+1} E_\alpha [(g(\eta_0)/{\eta_0}) \psi(\eta)] \\
&=& E_\alpha[(g(\eta_0)/{\eta_0}) \psi(\eta)] + 
((1-\lambda)^{R}-1)E_\alpha[(g(\eta_0)/{\eta_0}) \psi(\eta)]\\
&=&E_\alpha[\h \psi(\eta)] + J_3.
\end{eqnarray*}

Hence, we have that
\begin{equation}
\label{via}
E_\alpha[(\h - \l\phi)\psi] \ =\  -(I_2 +I_3 +I_4 +J_1+J_3) .
\end{equation}
To show the bound in (\ref{bounds1}), by the variational
characterization of $\|\cdot\|_{-1,\lambda}$ (cf. (\ref{var_expression})), we need only verify
$$|I_2 +I_3 +I_4 +J_1+ J_3| \ \leq \ C\|\psi\|_{1,\lambda}$$
for some constant $C=C(\alpha, p)$.

To this end, observe, by Schwarz inequality,
$$
|I_3| \ \leq \ \bigg (\sum_j |c|^2 p(j)\bigg )^{1/2} \bigg (\sum_j \alpha
E_\alpha[g(\eta_0)\frac{\eta_0-1}{\eta_0}(\psi(\eta^{0,-j}) -
\psi(\eta))^2]p(j)\bigg )^{1/2} \ \leq \ C\|\psi\|_1
$$
for a constant $C=C(\alpha)$ as $p(j)\leq 2s(-j)$.
Also, 
\begin{eqnarray*}
|I_4| &\leq& \bigg (\sum_j |c|^2 p(j)\bigg )^{1/2} \bigg (
\sum_j E_\alpha[g(\eta_0)\eta_0]E_\alpha[\frac{g(\eta_0)}{\eta_0}(\psi(\eta ) - \psi(\tau_{-j}( \eta^{0,-j})))^2]p(j)\bigg )^{1/2}\\
& \leq& 2|c| \sqrt{E_\alpha[g(\eta_0)\eta_0]} \|\psi\|_1^2.
\end{eqnarray*}

For the second term $I_2$, note, with Lemma \ref{computations},
\begin{eqnarray*}
&&E_\alpha[(g(\eta_i) -
{\eta_{i+j}}\frac{g(\eta_0)}{\eta_0})\psi(\eta)]\\
&&\ \ \ \ = E_\alpha[g(\eta_0)\frac{\eta_0-1}{\eta_0}(\psi(\eta^{0,i}) - \psi(\eta)))
+ (g(\eta_0)-g(\eta_0)\frac{\eta_{i+j} +1}{\eta_0})\psi(\eta)]\\
&&\ \ \ \ = E_\alpha[g(\eta_0)\frac{\eta_0-1}{\eta_0}(\psi(\eta^{0,i}) - \psi(\eta)))
+ g(\eta_0)(\psi(\eta)-\psi(\tau_{-(i+j)}(\eta^{0,-(i+j)})))].
\end{eqnarray*}
Then,
\begin{eqnarray*}
|I_2| &\leq&
\bigg (\sum_j\sum_{{|i| \leq R} \atop {i\neq 0,-j}}
|c|^2 p(j)\bigg )^{1/2}
\bigg (\sum_j\sum_{|i|\leq R} 2\alpha
E_\alpha[g(\eta_0)\frac{\eta_0-1}{\eta_0}(\psi(\eta^{0,i}) - \psi(\eta)))^2]p(j)\\
&&\ \ \ \ \ \ \ + 2E_\alpha[g(\eta_0)\eta_0]E_\alpha[\frac{g(\eta_0)}{\eta_0}(\psi(\eta) -
\psi(\tau_{-(i+j)}(\eta^{0,-(i+j)})))^2]p(j)\bigg )^{1/2}.
\end{eqnarray*}
Note, as $s$ is irreducible, $u\in \Z$ can be written $u = \sum_{k=1}^m l_k$ for
points $l_k$ in the support of $s$, $s(l_k)>0$.  Let $r_0=0$ and $r_k =
\sum_{n=1}^kl_n$ for $1\leq k\leq m$.  Then, with Lemma \ref{computations},
\begin{eqnarray*}
|E_\alpha[g(\eta_0)\frac{\eta_0-1}{\eta_0}(\psi(\eta^{0,u}) -
\psi(\eta)))^2]| &=& \alpha| E_\alpha[ (\psi(\eta +\delta_u) -\psi(\eta
+\delta_0))^2]|\\
&\leq& m\alpha\sum_{k=0}^{m-1} E_\alpha[
(\psi(\eta +\delta_{r_k})-\psi(\eta+\delta_{r_{k+1}}))^2]\\
&\leq &  C\|\psi\|_1^2
\end{eqnarray*}
for some constant $C=C(p,m)$ as $p$ is finite-range.
Also,
\begin{eqnarray*}
&&|E_\alpha[\frac{g(\eta_0)}{\eta_0}(\psi(\eta) -
\psi(\tau_{u}(\eta^{0,u})))^2]|\\
&&\ \ \ \leq m \sum_{k=0}^{m-1}
E_\alpha[\frac{g(\eta_0)}{\eta_0}(\psi(\tau_{r_k}(\eta^{0,r_k}))-\psi(\tau_{r_{k+1}}(\eta^{0,r_{k+1}})))^2]\\
&&\ \ \ = m\sum_{k=0}^{m-1}
E_\alpha[\frac{g(\eta_0)}{\eta_0}(\psi(\eta)-\psi(\tau_{l_{k+1}}(\eta^{0,l_{k+1}})))^2]
\ \leq \  C\|\psi\|_1^2
\end{eqnarray*}
for a $C=C(p,m)$ again as $p$ is finite-range.
Then, as the sums in $I_2$ are finite, $I_2$ is bounded
$|I_2|\leq C\|\psi\|_1^2$ for some constant $C=C(\alpha,p)$.

To bound $J_1$, we use the resolvent bound (\ref{resolventbound}).
Namely, as $\{g(\eta_i)-\alpha: i\in\Z^d\}$ is an orthogonal family,
\begin{eqnarray*}
&&|\sum_j\sum_{i\geq R+1}\a_{i+j,i} E_\alpha [(g(\eta_{i})-\alpha)\psi(\eta)]p(j)|\\
&&\ \ \ \ \ \leq
\|\sum_j\sum_{i\geq R+1}\a_{i+j,i} (g(\eta_{i})-\alpha)p(j)\|_{-1,\lambda} \|\psi\|_{1,\lambda}\\
&&\ \ \ \ \ \leq \bigg (\frac{1}{\lambda}\sum_j\sum_{i\geq R+1}\nabla^2 a_{i+j,i}E_\alpha[(g(\eta_1)-\alpha)^2]p(j)\bigg
)^{1/2} \|\psi\|_{1,\lambda}
\end{eqnarray*}
and, as $\{(\eta_{i+j}-\rho)\frac{g(\eta_0)}{\eta_0}: i\geq R+1\}$ is
also an orthogonal collection,
\begin{eqnarray*}
&&|\sum_j\sum_{i\geq
 R+1}\a_{i+j,i}E_\alpha[(\eta_{i+j}-\rho)\frac{g(\eta_0)}{\eta_0})\psi(\eta)]p(j)|\\
&&\ \ \ \ \ \ \ \ \ \ \leq \bigg (\frac{1}{\lambda}
\sum_j\sum_{i\geq
 R+1}\nabla^2 a_{i+j,i}E_\alpha[(\eta_{1}-\rho)^2\frac{g^2(\eta_0)}{\eta^2_0}]p(j)\bigg
 )^{1/2}\|\psi\|_{1,\lambda};
\end{eqnarray*}
then, $|J_1|\leq C(\lambda^{-1}\sum_{i\geq 1}\nabla^2
a_{i+1,i})\|\psi\|_{1,\lambda}$ for some $C=C(\alpha,p)$ as $i+j\geq
1$ for $i\geq R+1$ beyond the range of $p$.

Finally, $J_3$ is bounded by the resolvent bound (\ref{resolventbound})
$$|J_3| \ \leq \ \bigg[\frac{|(1-\lambda)^R
  -1|}{\lambda}\|g(\eta_0)/\eta_0\|^2_0\bigg]^{1/2} \|\psi\|_{1,\lambda}.$$

Putting these estimates together, using a form of Schwarz--relation $2ab =
\inf_{\epsilon} \epsilon^{-1}a^2 +\epsilon b^2$--we obtain, 
for a constant $C=C(\alpha,p)$,
$$|I_2 +I_3+I_4 + J_1+J_3| \ \leq \ C\bigg(1+ \
\frac{1}{\lambda}\sum_{i\geq 1} \nabla^2 a_{i+1,i} 
\bigg)^{1/2}\|\psi\|_{1,\lambda}.
$$

Now, by direct computation, we have that
$$\sum_{i\geq 1} a_i^2 \ = \ \frac{|c|^2}{\lambda} \ \ \ 
{\rm \ and \ \ \ \ } \frac{1}{\lambda}\sum_{i\geq
  1}\nabla^2 a_{i+1,i} \ = \ \frac{\lambda^2
  |c|^2}{\lambda}\frac{1}{\lambda (2-\lambda)}$$
which shows (\ref{bounds1}) via (\ref{via}).

To show (\ref{bounds2}), we observe
$
\lambda
\|\phi\|^2_{L^2}  =  \lambda (|c|^2/\lambda) = |c|^2$
and
\begin{eqnarray*}
\|\phi\|_{1}^2 & = & \alpha\sum_j \sum_i \nabla^2 a_{i+j,i}s(j)\\
&&\ \ \ + \sum_j E_\alpha[\frac{g(\eta_0)}{\eta_0}\bigg (\sum_i
\a_{i,i+j}(\eta_{i+j}-\rho)+ \sum_j\a_{0,-j}\bigg )^2]p(j)\\ 
&\leq & C + C\sum_j \sum_{i\geq 1} \nabla^2 a_{i+j,i} \ \leq \  C'
\end{eqnarray*}
for some constants $C=C(\alpha,p)$ and $C'=C'(\alpha,p)$ using, as before, the orthogonality of
$\{({g(\eta_0)}/{\eta_0})(\eta_{i+j}-\rho)\}$.

\vskip .3cm
{\bf 2.1.2 Estimates in $d\geq 3$ under (SP).}  From (\ref{lphipsi}), for the function
$\phi(\eta)= (\eta_{j_0}-\rho)/(\rho p(j_0))$ and $\psi\in \D$, we have
$\rho p(j_0)E_\alpha[(\l\phi)\psi]$ equals
\begin{eqnarray*}
&& -E_\alpha[(g(\eta_{j_0}) -
\eta_{2j_0}g(\eta_0)/\eta_0)\psi]p(j_0) + E_\alpha[(g(\eta_{2j_0})- \eta_{j_0}g(\eta_0)/\eta_0)\psi]p(-j_0)\\
&&\ \ +\sum_{j\neq \pm j_0}
\bigg\{E_\alpha[(g(\eta_{j_0-j})-\eta_{j_0}g(\eta_0)/\eta_0)\psi]
-E_\alpha[(g(\eta_{j_0})-\eta_{j_0+j}g(\eta_0)/\eta_0)\psi]\bigg\}p(j)\\
&&\ \ - \alpha a_{j_0}E_\alpha[\psi(\eta+\delta_{j_0}) - \psi(\eta +
\delta_0)]p(-j_0)\\
&&\ \  + a_{j_0}E_\alpha[g(\eta_0)(\psi(\eta)-\psi(\tau_{-j_0}(\eta^{0,-j_0})))]p(j_0).\end{eqnarray*}
As we can take $E_\alpha[\psi]=0$ without loss of generality, with
Lemma \ref{computations}, 
$\rho p(j_0)E_\alpha[(\l\phi)\psi]$ equals
\begin{eqnarray*}
&& \rho E_\alpha[(g(\eta_0)/\eta_0)\psi]p(j_0)
-E_\alpha[(g(\eta_{j_0}) -\alpha)\psi]p(j_0)
+E_\alpha[(\eta_{2j_0}-\rho)(g(\eta_0)/\eta_0)\psi]p(j_0)\\
&&\ \ + \sum_{j\neq \pm j_0}
\bigg\{E_\alpha[(g(\eta_{j_0-j})-\eta_{j_0}g(\eta_0)/\eta_0\psi]
-E_\alpha[(g(\eta_{j_0})-\eta_{j_0+j}g(\eta_0)/\eta_0\psi]\bigg\}p(j)\\
&&\ \
-a_{j_0}E_\alpha[g(\eta_0)\frac{\eta_0-1}{\eta_0}\psi(\eta^{0,j_0}) - 
\psi(\eta)]p(-j_0)\\
&&\ \  + a_{j_0}E_\alpha[g(\eta_0)(\psi(\eta)-\psi(\tau_{-j_0}(\eta^{0,-j_0})))]p(j_0) \\
&&= \rho p(j_0)E_\alpha[(g(\eta_0)/\eta_0)\psi] +K_1 +K_2+K_3 +K_4 +K_5.
\end{eqnarray*}

Hence, to show $\|\h-\l\phi\|_{-1,\lambda}<\infty$, by the variational
characterization (cf. (\ref{var_expression})), we need only show that
$$|E_\alpha[(\h-\l \phi)\psi]| \ = \ (\rho p(j_0))^{-1}|K_1 +K_2 +K_3 +K_4+K_5| \ \leq \
C\|\psi\|_1$$
for some constant $C=C(\alpha,p)$.
To this end, the terms $K_3, K_4$ and $K_5$ are handled analogously as
$I_2, I_3$ and $I_4$ above in
the $d=1$ case.  To bound $K_1$ and $K_2$, we invoke the following
result.

\begin{prop} 
\label{SXprop}
Consider $d\geq 3$ reference frame processes such that $g$
  satisfies assumption (SP).  Let $f$ be a $L^4(Q_\alpha)$ function supported on a finite
  number of vertices of $\Z^d\setminus\{0\}$ which is mean-zero,
  $E_\alpha[f]=0$.  Then, $\|f\|_{-1} < \infty$.
\end{prop}

{\it Proof.}  The proof is virtually the same as for Theorem 1.2 \cite{SX}.  One can
bound in terms of a constant $C=C(f,\alpha,p,d)$ that
$$E_\alpha[f\psi] \ \leq \ 
C \bigg [\sum_j \sum_{i\neq 0}
E_\alpha[g(\eta_i)(\psi(\eta^{i,i+j})-\psi(\eta))^2]s(j)\bigg
]^{1/2}\ \leq \ C\|\psi\|_1$$
by straightforwardly avoiding the orign.  This relation gives the
desired statement. \qed

Note now Proposition \ref{SXprop} directly applies to $K_1$.  For $K_2$, we first condition on
$\eta_0$ to get 
\begin{eqnarray*}
|K_2| &=&
p(j_0)|E_\alpha[E_\alpha[(\eta_{2j_0}-\rho)\psi(\cdot ;\eta_0)|\eta_0](g(\eta_0)/\eta_0)]|\\
&\leq& p(j_0) E_\alpha[C\|\psi(\cdot;\eta_0)\|_1 (g(\eta_0)/\eta_0)]\\
&\leq& Cp(j_0) E_\alpha[(g(\eta_0)/\eta_0)^2]^{1/2} \|\psi\|_1
\end{eqnarray*}
where $\psi(\cdot;\eta_0)$ denotes $\psi$ as a function of 
$\{\eta_i: i\neq 0\}$ with $\eta_0$ fixed.
This finishes the proof of Theorem \ref{thm1} in this case. \qed

\sect{Proof of Theorem \ref{thm2}}

We first define
the notion of a positively associated stationary increments $L^2$ process
$N(t)$.  
This
is an $L^2$ process where
$$E[\phi(N(t+s)-N(t))\psi (N(t))] \geq E[\phi(N(s))]E[\psi(N(t))]$$
for all $\phi$ and $\psi$ increasing.
For such processes we have the Newman-Wright result (cf. \cite{N}).
\begin{theorem}
\label{NW}
Suppose $N(t)$ is an $L^2$ process with 
positively associated stationary increments such
that the limit exists
$$lim_{t\rightarrow \infty} \frac{1}{t}E\bigg [(N(t)-E[N(t)])^2\bigg ] \ = \ \sigma^2 <\infty.$$
Then, we have weak convergence to Brownian motion in Skorohod space, 
$$ \frac{1}{\sqrt{\lambda}}\bigg(N(t) -E[N(t)]\bigg) \ \rightarrow \
\sigma B(t).$$
\end{theorem}

The strategy will now be to verify that the tagged position
$x(t)$ in $d=1$ under the assumptions of Theorem \ref{thm2} has associated increments and that its variance scales
diffusively so that the Newman-Wright statement applies.

The following is a useful coupling which essentially says
adding more particles to the system
only slows down the tagged particle.
\begin{lemma}
\label{coupling}
Under the assumptions on $g$ and $p$ in dimension $d=1$ in Theorem
\ref{thm2},
we can couple two
two copies of the joint 
process, $(x^1(t),\xi^1(t))$ and $(x^2(t),\xi^2(t))$ where
$\xi^1(0) \leq \xi^2(0)$ 
coordinatewise and $x^1(0) \geq x^2(0)$
so that at all later time $t$, $x^1(t)\geq x^2(t)$.
\end{lemma}

{\it Proof.}
We make the coupling so that when
an $\xi^1$ particle moves, a corresponding $\xi^2$ particle also moves
to the right, and also when $x^2$ would move ahead of $x^1$ then $x^1$
also moves.

More carefully, at vertex $x\neq x^1, x^2$, the basic coupling
applies--with rate $g(\xi_x^1)$ a particle from $x$ in both systems
moves; and with rate $g(\xi_x^2)-g(\xi_x^1)$ a particle from $x$ in
system $2$ moves.

When $x^1\neq x^2$, with rate
$g(\xi^1_{x^1})(\xi^1_{x^1}-1)/\xi^1_{x^1}$ a non-tagged particle
in system $1$ and a particle in system $2$ moves from location $x^1$;
with rate $g(\xi^1_{x^1})/\xi^1_{x^1}$ the
tagged particle from system $1$ and a particle from system $2$ at
$x^1$ moves;
and with rate $g(\xi^2_{x^1})-g(\xi^1_{x^1})$ a particle in system $2$
at $x^1$
moves.  

With respect to location $x^2$, with rate
$g(\xi^1_{x^2})(\xi^1_{x^2}-1)/\xi^1_{x^2}$ a particle from system $1$
and a non-tagged particle in system $2$ moves from $x^2$; with rate
$g(\xi^2_{x^2})/\xi^2_{x^2}$ the tagged particle in system $2$ and a
particle in system $1$ move from location $x^2$; with rate
$g(\xi^1_{x^2})/\xi^1_{x^2} - g(\xi^2_{x^2})/\xi^2_{x^2}$ a particle
in system $1$ moves from $x^2$; with rate
$g(\xi^2_{x^2})(\xi^2_{x^2}-1)/\xi^2_{x^2} -
g(\xi^1_{x^2})(\xi^1_{x^2}-1)/\xi^1_{x^2}$ a non-tagged particle moves
from system $2$ at $x^2$.

When $x^1=x^2=x$, with rate $g(\xi^1_x)(\xi^1_x -1)/\xi^1_x$ a
non-tagged particle from $x$ in both systems moves; with rate
$g(\xi^2_x)/\xi^2_x$ both tagged particles move; with rate
$g(\xi^1_x)/\xi^1_x - g(\xi^2_x)/\xi^2_x$ the tagged particle in
system $1$ and a non-tagged particle in system $2$ moves; with (the
remaining) rate
$$g(\xi^2_x)\frac{\xi^2_x -1}{\xi^2_x} - g(\xi^1_x)\frac{\xi^1_x
  -1}{\xi^1_x} - \frac{g(\xi^1_x)}{\xi^1_x} +\frac{g(\xi^2_x)}{\xi^2_x} \ = \ g(\xi^2_x) - g(\xi^1_x)$$
a non-tagged particle in system $2$ moves.

We omit the generator formulation. \qed

The next lemma owes some intuition to Theorem 2 \cite{Kipnis}.
\begin{lemma}
\label{posassoc}
In $d=1$, under the assumptions of Theorem \ref{thm1}, the $L^2$
process $x(t)$ under equilibrium $Q_\alpha$ has positively associated stationary increments.
\end{lemma}

{\it Proof.} From (\ref{v_decomp}), clearly
$x(t)$ is an $L^2$ process.  Also under equilibrium $Q_\alpha$,
$x(t)$ has stationary increments.  Consider now the sequence, for
increasing $\phi$ and $\psi$,
\begin{eqnarray*}
E_\alpha [\phi(x(t+s)-x(t))\psi(x(t))] 
&=& E_\alpha [ \psi(x(t))E_{\eta(t)}[\phi(x(s))]]\\
&=& E^*_\alpha[\psi(x^*(0)-x^*(t))E_{\eta^*(0)}[\phi(x(s))]]\\
&=& E^*_\alpha[E^*_{\eta^*(0)}[\psi(x^*(0)-x^*(t))]E_{\eta^*(0)}[\phi(x(s))]]\\
&=& \int E^*_{\eta^*}[\psi(x^*(0)-x^*(t))]E_{\eta^*}[\phi(x(s))]dQ_\alpha (\eta^*)
\end{eqnarray*}
where in the second step we note $x(0)=0$, and reverse time at $t$ with
$x^*(u) = x(t-u)$, $\eta^*(u)= \eta(t-u)$,
and $E^*_\alpha$ and $E^*_{\eta^*}$ denotes expectation with respect to the reversed
process with initial distribution $Q_\alpha$ and state $\eta^*$
respectively.

Consider the functions $E^*_{\eta^*}[\psi(x^*(0)-x^*(t))]$,
and $E_{\eta^*}[\phi(x(s))]$ as functions of $\eta^*$.  Both are decreasing
coordinatewise by the coupling in Lemma \ref{coupling}.  
Indeed, from the coupling, we see that, by increasing $\eta^*$ by one particle,
$x^*(0)-x^*(t)$ decreases (recall that the reversed $*$ process moves
to the left), and $x(t)$ decreases.
In other words, both functions decrease coordinatewise in $\eta^*$.

With this monotonicity, the associated property follows from
the standard FKG inequality for product measures (see Liggett \cite{Liggett}).
\qed

We now turn to an analysis of the variance.  Recall the definition of
$\f$ (cf. near (\ref{martdecomp})).
\begin{lemma}
\label{varlemma}
In all $d\geq 1$, 
the variance $V(t)= E_\alpha[|x(t)-E_\alpha[x(t)]|^2]$ satisfies
$$V(t) \ = \ \frac{\alpha}{\rho}\sum |j|^2p(j)t + 2\int_0^t
E_\alpha[x(s)\cdot \f(\eta(s))]ds.$$
\end{lemma}

{\it Proof.}
We continue the sequence (\ref{martsequence}).  Write
\begin{eqnarray*}
V(t) &=& \frac{\alpha}{\rho}\sum |j|^2p(j)t + 2\int_0^t
E_\alpha [M(s)\cdot \f(\eta(s))]ds + E_\alpha[|A(t)|^2]\\
&=& \frac{\alpha}{\rho}\sum |j|^2p(j)t + 2\int_0^t E_\alpha[x(s)\cdot \f(\eta(s))]ds
-2E_\alpha[A(s)\cdot \f(\eta(s))]
+ E_\alpha[|A(t)|^2]\\
&=& \frac{\alpha}{\rho}\sum |j|^2p(j)t + 2\int_0^t E_\alpha[x(s)\cdot \f(\eta(s))]ds.
\end{eqnarray*}
Here, in the second line, we use $x(t) - E_\alpha[x(t)] = M(t) + A(t)$, and
in the last line that $|A(t)|^2 = 2\int_0^t
\int_0^s\f(\eta(r))\cdot \f(\eta(s))drds = 2\int_0^t A(s)\cdot \f(\eta(s))ds$. \qed

\begin{lemma}
\label{superadd}
In $d=1$, the variance
$V(t)= E_\alpha[(x(t)-E_\alpha[x(t)])^2]$ is super-additive, and so
the limit $\lim_{t\rightarrow \infty} V(t)/t$ exists.
\end{lemma}

{\it Proof.}
We study the term $E_\alpha [x(s)\f(\eta(s))]$ appearing in the variance
expression in Lemma \ref{varlemma}.  Reverse time at $s$ (using the notation
given in proof of Lemma \ref{posassoc}), to
obtain
\begin{eqnarray}
\label{U}
E_\alpha [x(s)\f(\eta(s))] &=& E_\alpha^*[(x^*(0) - x^*(s))\f(\eta^*(0))]
\end{eqnarray}

Now, write, using (\ref{U}) and the variance decomposition in Lemma \ref{varlemma}, that
\begin{eqnarray*}
V(t+r)-V(r) - V(t) &=& 2\int_0^{t}
E_\alpha^*[(x^*(s)-x^*(s+r))\f(\eta^*(0))]ds\\
&=& 2 \int_0^t E_\alpha^* [ E^*_{\eta^*(s)}[x^*(0)-x^*(r)]
\f(\eta^*(0))]ds\\
&=& 2\int_0^t E_\alpha [ E^*_{\eta(0)}[x^*(0)-x^*(r)] \f(\eta(s))]ds\\
&=& 2 E_\alpha [ E^*_{\eta(0)}[x^*(0)-x^*(r)]
E_{\eta(0)}[A(t)] ]]\\
&=& 2 E_\alpha [ E^*_{\eta(0)}[x^*(0)-x^*(r)]
E_{\eta(0)}[x(t)-E_\alpha [x(t)] ]].
\end{eqnarray*}
Here, we shifted variables $s$ to $s+r$ in the first line, conditioned
up to time $s$ in the second line, reversed time at $s$ in the third,
and used the martingale decomposition $x(t)-E_\alpha [x(t)]=M(t)+A(t)$
in the last
line.

The last product, as in proof of Lemma \ref{posassoc}, is the product of
decreasing functions of $\eta$.  Now use FKG inequality to finish the
proof. \qed

We are now ready to prove Theorem \ref{thm2}.
\vskip .2cm

{\it Proof of Theorem \ref{thm2}.} 
Under the assumptions of Theorem \ref{thm2}, we can invoke the
Newman-Wright principle:  By Lemma \ref{posassoc}, $x(t)$ has
positively associated increments.  By Lemma \ref{superadd}, the limit
$\lim_{t\rightarrow \infty}
V(t)/t =\sup_{t\geq 1} V(t)/t$ exists; and, by Theorem
\ref{thm1}, $\sup_{t\geq 1}
V(t)/t <\infty$. 

Finally, to show the limit $\lim_{t\rightarrow \infty} V(t)/t>
\alpha/\rho$ we need only show by superadditivity, noting (\ref{U}), that
\begin{equation}
\label{helpingU}
\int_0^1E_\alpha[x(s)\f(\eta(0))]ds \ >\ 0.
\end{equation}
But, one can write
$$
E_\alpha[x(s)\f(\eta(s))] \ =\  E_\alpha\bigg[(x(t)
-E_\alpha[x(t)])\frac{g(\eta_0)}{\eta_0}\bigg]
\ =\ \frac{\alpha}{\rho}\bigg\{E'_\alpha[x(t)] - E_\alpha[x(t)]\bigg\}
$$
after some algebra where $E'_\alpha$ is the process expectation with
respect to initial distribution $Q'_\alpha = \prod_{i\neq
   0}\mu_\alpha \times \mu'_\alpha$ and
 $$\mu'_\alpha (k)
 \ = \  \frac{1}{Z'_\alpha} \frac{\alpha^{k-1}}{g(1)\cdots g(k-1)} \ \ \
   {\rm for \ } k\geq 1$$
 with normalization
 $Z'_\alpha$.  The interpretation is that $\mu'_\alpha$ puts a particle
 at the origin and distributes other particles there according to
 $\mu_\alpha$.
 It is a straightforward computation, under Assumption (ID), that
 $\mu'_\alpha << \mu^0_\alpha$ in stochastic order, and so can couple
  two joint systems starting from $Q_\alpha$ and $Q'_\alpha$ so that the
 tagged particle under $Q_\alpha$ always is ahead of its counterpart
 under $Q'_\alpha$. 

Since the inequality $\mu'_\alpha <<\mu_\alpha$ is strict, with
positive probability, the $Q'_\alpha$ system has strictly less
particles than the $Q_\alpha$ system initially at the origin.
It is not hard now to construct a
situation with positive probability, as all clocks are exponential, 
where the tagged particle
positions differ at some time $0<t\leq 1/2$, and that this difference
is maintained up to time $t=1$.  Hence, (\ref{helpingU}) holds.
\qed

\vskip .2cm
Last, we discuss briefly here a ``particle-level'' approach for why 
$\sup_{t\geq 1} V(t)/t <\infty$ in the case of Theorem \ref{thm2} should be
true.
One needs only bound the
``drift'' part $2\int_0^t
E_\alpha^*[(x^*(0)-x^*(s))\f(\eta^*(s))]ds$, or show
$$E_\alpha[(x(t)-E_\alpha[x(t)])\frac{g(\eta_0(0))}{\eta_0(0)}] \ = \ 
\frac{\alpha}{\rho}\big \{ E'_\alpha[x(t)] - E_\alpha[x(t)]\big\} \ =
\ O(1).$$
To bound uniformly the difference between
the tagged particles in expected value, 
the key point would be to handle
the influence of ``extra'' particles at the origin in
the $Q_\alpha$ system which could ``slow down'' the $Q_\alpha$ tagged
particle and make a large difference.
However, these extras, though not quite ``second-class'' particles,
should in the long term behave like them and move much slower
than a tagged particle, and so their influence should be negligible in the
limit.  Making this precise technically however seems difficult.  We remark though
that these ideas led in part to other interesting ``point-of-view
shifts'' problems \cite{HL}.

\vskip .2cm
{\bf Acknowledgement.}  I would like to thank P. Ferrari, N. Jain,
C. Landim and C. Lee for useful discussions.

 \bibliographystyle{plain}

\vskip .5cm
Sunder Sethuraman

400 Carver Hall

Iowa State University

Ames, IA \ 50011

sethuram@iastate.edu
\end{document}